\documentclass[leqno,11pt, a4]{amsart}
\usepackage[utf8]{inputenc}
\usepackage{amsthm}
\usepackage{amsmath}
\usepackage{mathtools}
\usepackage{enumitem}
\usepackage{amsfonts}
\usepackage{amssymb}
\usepackage{graphicx}
\usepackage{pdfsync}
\def\R{\text{$\mathbb{R}$}}

\def\d1#1#2{\frac{d#1}{d#2}}

\def\p1#1#2{\frac{\partial #1}{\partial #2}}

\def\C{\text{$\mathbb{C}$}}
\def\N{\text{$\mathbb{N}$}}
\def\Z{\text{$\mathbb{Z}$}}

\newcommand{\di}{\text{d}}

\newcommand{\sph}{\mathbb{S}}

\newcommand{\g}[1]{\mathfrak{#1}}

\newcommand {\lag}{\mathcal{L}}

\newcommand{\dcp}[2]{\left\langle #1|#2\right\rangle}
\usepackage{comment}

\newcommand {\der}[2]{\frac{d#1}{d#2}}
\newcommand {\im}[1]{\mathfrak{\mathfrak{I}m}(#1)}
\newcommand{\Ad}{\text{Ad}}

\title[]{Complex Structures on Product Manifolds}
\author{Leonardo Biliotti and Alessandro Minuzzo }
\address{Dipartimento di Scienze Matematiche, Fisiche e Informatiche \\
          Universit\`a di Parma (Italy)}
\email{leonardo.biliotti@unipr.it}
\email{alessandro.minuzzo@unipr.it}
\thanks{The first author was partially supported by PRIN  2017
   ``Real and Complex Manifolds: Topology, Geometry and holomorphic dynamics ''
   and ``National Group for Algebraic and Geometric Structures, and their Applications'' (GNSAGA - INDAM)}
\begin{document}
\newtheorem{thm}{Theorem}[section]
\newtheorem*{red}{Reduction Principle for Polar Actions}
\newtheorem*{redc}{Reduction Principle for Coisotropic actions}
\newtheorem*{redalmost}{Reduction Principle for Almost Homogeneous actions}
\newtheorem*{eqt}{Equivalence Theorem for Hamiltonian actions}
\newtheorem*{equivalence-nonHamiltonian}{Equivalence Theorem}
\newtheorem*{cauzzi}{Restriction Lemma}
\newtheorem{prop}[thm]{Proposition}
\newtheorem{lemma}[thm]{Lemma}
\newtheorem{cor}[thm]{Corollary}
\theoremstyle{definition}
\newtheorem{defini}[thm]{Definition}
\newtheorem{notation}[thm]{Notation}
\newtheorem{exe}[thm]{Example}
\newtheorem{conj}[thm]{Conjecture}
\newtheorem{prob}[thm]{Problem}
\theoremstyle{remark}
\newtheorem{rem}[thm]{Remark}
\newcommand{\thistheoremname}{}
\theoremstyle{plain}
\newtheorem*{genericthm*}{\thistheoremname}
\newenvironment{namedthm*}[1]
  {\renewcommand{\thistheoremname}{#1}%
   \begin{genericthm*}}
  {\end{genericthm*}}
\newcommand{\restr}[1]{\vert_{#1}}
\newcommand{\act}[2]{#1 \circlearrowright #2}
\newcommand{\actg}{G \circlearrowright X}
\newcommand{\actgl}{G \circlearrowright M}
\newcommand{\orbg}{G(p)}
\newcommand{\euclideo}{\langle \cdot , \cdot \rangle}
\def\dq{\mathrm{d}_{M/G}}
\def\dgr{\mathrm{d}_{M^H/N(H)}}
\newcommand{\princ}{M_{\mathrm{princ}}}
\begin{abstract}
Let $M_i$, for $i=1,2$, be a Kähler manifold, and let $G$ be a Lie group acting on $M_i$ by Kähler isometries. Suppose that the action admits a momentum map $\mu_i$ and let $N_i:=\mu_i^{-1}(0)$ be a regular level set. When the action of $G$ on $N_i$ is proper and free, the Meyer--Marsden--Weinstein quotient $P_i:=N_i/G$ is a Kähler manifold and $\pi_i:N_i\to P_i$ is a principal fiber bundle with base $P_i$ and characteristic fiber $G$. In this paper, we define an almost complex structure for the manifold $N_1\times N_2$ and give necessary and sufficient conditions for its integrability. In the integrable case, we find explicit holomorphic charts for $N_1\times N_2$. As applications, we consider a non integrable almost-complex structure on the product of two complex Stiefel manifolds and the infinite Calabi-Eckmann manifolds $\sph^{2n+1}\times S(\mathcal{H})$, for $n\geq 1$, where $S(\mathcal{H})$ denotes the unit sphere of an infinite dimensional Hilbert space $\mathcal{H}$.
\end{abstract}
\maketitle
\section{Introduction}
In a classical paper by E. Calabi and B. Eckmann \cite{calheck} it is shown that the compact manifold $\sph^{2n+1}\times \sph^{2m+1}$ admits a complex structure coming from the natural $U(1)$-action on $\C^{n}\times \C^m$, but does not admit a Kähler structure for $n\geq 1$, as an immediate cohomological computation shows. For example, the second de Rham cohomology group $H^2_{dR}(\sph^{2n+1}\times\sph^{2m+1})=\{0\}$, precisely when $n\geq 1$ or $m\geq 1$.

We generalize the idea of their construction by considering two Kähler manifolds $(M_i,h_i,J_i,\omega_i)$, $i=1,2$, where $h_i$ is the Riemannian metric, $J_i$ the complex structure, $\omega_i$ the compatible symplectic structure and letting $G$ be a Lie group acting on each $M_i$ by \textit{Kähler isometries}. This means that for all $g\in G$ we have 

$$g^*h_i=h_i, \qquad g_*\circ J_i=J_i\circ g_*,\qquad g^*\omega_i=\omega_i\ ,$$
where by $g^*$ and $g_*$ we denote respectively the pull-back and the push-farward, associated to the diffeomorphism induced by $g$. 
We suppose that each action admits a \textit{momentum map} $\mu_i:M_i\to \g{g}^*$, where $\g{g}:=\text{Lie}(G)$, and let $N_i$ be a regular level set of it. We further suppose that we can reduce $N_i$ to its orbit space $P_i$ by taking the quotients by the isotropy subgroup $F_i$ of the momentum, relative to the co-adjoint action of $G$ on $\g{g}^*$. This is the Meyer--Marsden--Weinstein quotient, and is known to be a Kähler manifold $(P_i,\tilde h_i,\tilde J_i,\tilde\omega_i)$, where ${h_i|}_{N_i}=\pi_i^*\tilde h_i$, ${\omega_i|}_{N_i}=\pi_i^*\tilde\omega_i$ and $\pi_i\circ J_i=\tilde J_i\circ \pi_i$. Moreover, the natural projection $\pi_i:N_i\to P_i$ is a principal fiber bundle with base $P_i$ and characteristic fiber $F_i$. For the definition of the momentum map and its use in the Meyer--Marsden--Weinstein reduction, we refer to \cite{Ca}.

There is a natural decomposition of the tangent space of $N_i$ as we now describe. For further details see \cite{KN}. The projection defines a distribution inside $TN_i$, the so called \textit{vertical distribution}, given by the kernel of its differential. At each point the vertical distribution is  isomorphic (not canonically) to the Lie algebra of the isotropy subgroup, $\g{f_i}:=\text{Lie}(F_i)$, because we are in the hypothesis that the action of $F_i$ on $N_i$ is free. We denote the complementary subspace at $p_i\in N_i$ with respect to $T_{p_i}N_i$ as $H_{p_i}$, i.e. $T_{p_i}N_i\simeq H_{p_i}\oplus \g{f_i}$. We can therefore define a distribution $H_i:=\cup_{p_i\in N_i}H_{p_i}\subset TN_i$, the \textit{horiziontal distribution}, and we also have that for all $g\in G$ that $g_*H_i=H_i$, because we are considering an isometric action. This tells us that $H_i\subset TN_i$ is a \textit{principal connection} for the principal fiber bundle $\pi_i:N_i\to P_i$. 

In the following, we will consider the case in which $N_i:=\mu_i^{-1}(0)$ is a regular level set of the momentum map $\mu_i$, so that $F_1\simeq G\simeq F_2$. Moreover, we will use the splitting 

\begin{equation}\label{splitting}
T_{(p_1,p_2)}(N_1\times N_2)\simeq (H_{p_1}\oplus \g{g})\oplus (H_{p_2}\oplus \g{g})\ .
\end{equation}

\begin{defini}\label{almostcomlex}
    
We define the following endomorphism $J$ of $T(N_1\times N_2)$ as 

$$J:=(J_1\oplus T_1)\oplus (J_2\oplus T_2)\ ,$$

\noindent where we have used the splitting (\ref{splitting}), and we have set for $\xi\in\g{g}$

$$T_1(\xi_{N_1}):=\xi_{N_2}\qquad T_2(\xi_{N_2}):=-\xi_{N_1}\ , $$
whith $\xi_{N_i}$ the infinitesimal generator of the action of $G$ on $N_i$, relative to $\xi$.
\end{defini}

Our first result is the following.

\begin{thm}\label{thm1}
The endomorphism $J$ defines an \textit{almost complex structure} on $N_1\times N_2$, i.e. $J^2=-\text{id}_{T(N_1\times N_2)}$. Furthermore, $J$ is integrable, i.e. complex, if and only if $G$ is abelian.
\end{thm}
\noindent By the celebrated Newlander--Nirenberg Theorem \cite{NN}, the integrability of the almost-complex structure $J$ on $N_1\times N_2$ is equivalent to the vanishing of the Nijenhuis tensor associated to $J$. This tensor is defined as 
$$N_J(X,Y):=[JX,JY]-[X,Y]-J[JX,Y]-[X,JY]\ ,$$
where $X,Y$ are smooth vector fields on $N_1\times N_2$.\\

We then study the action by the abelian group $(\C,+)$ defined as \begin{equation}\label{psi1}
    \Psi:\C\times (N_1\times N_2)\to (N_1\times N_2),\qquad (a+ib, (x,y))\mapsto (e^{2\pi i (a+b)}.x,e^{2\pi i (a-b)}.y)\ , 
    \end{equation}
    where the dot denotes the action on the two factors. We prove the following proposition, telling us that the map $\Psi$ is holomorphic.

\begin{prop}\label{psihol1}

\begin{enumerate}
    \item The map $\Psi_{(x,y)}:\C\to (N_1\times N_2)$, $\Psi_{(x,y)}(a+ib):=\Psi(a+ib,(x,y))$ is holomorphic with respect to the complex structure of $\C$ and the complex structure $J$ of definition \ref{almostcomlex};

    \item The map $\Psi_{a+ib}:(N_1\times N_2)\to (N_1\times N_2)$, $\Psi_{a+ib}(x,y):=\Psi(a+ib,(x,y))$ is holomorphic with respect to the complex structure $J$ of definition \ref{almostcomlex}.
\end{enumerate}

\end{prop}

\noindent Let now $\Lambda:=\{a+ib\ :\ a,b\in\mathbb{Z}/2\}$ be the lattice of $\C$ with cell $\left[0,\frac{1}{2}\right]\times\left[0,\frac{i}{2}\right] $. Then $\C/\Lambda\simeq \mathbb{T}^2=\sph^{1}\times \sph^1$ is the standard one dimensional complex torus. This is a compact abelian group and the quotient action $\Psi$ by $\mathbb{T}^2$ is free. Without using theorem \ref{thm1}, we are able to prove the following 
\begin{thm}\label{thm2}
    
Given the action (\ref{psi1}), the natural projection $\pi:(N_1\times N_2)\to (P_1\times P_2)$ is a principal holomorphic bundle with characteristic fiber $\mathbb{T}^2$. Moreover, the construction can be easily generalized to the case of an action of $\C^n/\Lambda_n=\mathbb{T}^{2n}$.
\end{thm}

\section{Proofs of the results}

\begin{proof} of \textbf{Theorem \ref{thm1}.} We start by proving that $J$ is an almost complex structure. Consider a basis element for $T_{(p_1,p_2)}(N_1\times N_2)\simeq (H_{p_1}\oplus \g{g})\bigoplus (H_{p_2}\oplus \g{g}) $ denoted as $(v_1, \xi_1,v_2, \xi_2)$, where $v_i\in H_{p_i}$ and $\xi_i\in\g{g}$ are nonzero elements. By definition \ref{almostcomlex}, and since $J_i$ is a complex structure on $H_{p_i}$, we have that
$$J^2(v_1, \xi_1,v_2, \xi_2)=J(J_1 v_1,\xi_2,  J_2 v_2, -\xi_1)=(J_1^2 v_1,-\xi_1,J_2^2 v_2,-\xi_2)=-(v_1,\xi_1,v_2,\xi_2) .$$

To prove integrability, we show the vanishing of the Nijenhuis tensor of $J$, denoted as $N_J$. We do this for each case comprised in the splitting (\ref{splitting}).
If $X_i$ is a horizontal vector field the integrability follows directly from the integrability of the $J_i$. Let $Y_1:=\xi_{N_1}$, $Y_2:=\eta_{N_2}$ be vertical vector fields, we have  
$$N_J(Y_1, Y_2)=[JY_1,JY_2]-[Y_1,Y_2]-J[JY_1,Y_2]-J[Y_1,JY_2]=-J[JY_1,Y_2]-J[Y_1,JY_2]=$$
$$=-J([\xi_{N_2},\eta_{N_2}]-[\xi_{N_1},\eta_{N_1}])=-J([\xi,\eta]_{N_1}-[\xi,\eta]_{N_2})\ .$$
Therefore, $N_J(Y_1, Y_2)=0$ if and only if $G$ is abelian.
By skew-symmetry and linearity we are left to check that $N_J(X_1,Y_i)=0$.

$$N_J(X_1, Y_1)=-[X_1,Y_1]-J[JX_1,Y_1]=-[X_1,Y_1]-J^2[X_1,Y_1]=[X_1,Y_1]-[X_1,Y_1]=0\ . $$
Moreover,

$$N_J(X_1, Y_2)=[JX_1,JY_2]-J[X_1,JY_2]=J[X_1,JY_2]-J[X_1,JY_2]=0\ .$$
In the last two equations we used the fact that for $X_i$ horizontal and $Y_i=\xi_{N_i}$ vertical we have 
$$[JX_i, Y_i]=[J_i X_i,Y_i]=J_i[X_i,Y_i]=J[X_i,Y_i]\ .$$
Indeed $J_i\circ g_*= g_*\circ J_i$ is equivalent to $\lag_{\xi_{N_i}}J_i=0$ for all $\xi\in\g{g}$, and hence $$[Y_i,J_iX_i]= \lag_{Y_i} (J_iX_i)=(\lag_{Y_i}J_i)(X_i)+J_i(\lag_{Y_i}X_i)=J_i(\lag_{Y_i}X_i)=J_i[Y_i,X_i]\ .$$  \end{proof}

\begin{proof} of \textbf{Proposition \ref{psihol1}.}

\begin{enumerate}
\item Consider the curve $t\mapsto (a+t+ib)\in \C$, so that $\der{}{t}(a+t+ib)|_{t=0}=1\in T_{(a+ib)}\C$. Therefore, 
$${(\di\Psi_{(x,y)})}_{a+ib}(1)=\der{}{t}\Psi_{(x,y)}(a+t+ib)\big|_{t=0}=(1_{N_1},1_{N_2})\ .$$
For the tangent vector $i$ consider the curve $t\mapsto (a+i(t+b))\in \C$, so that $\der{}{t}(a+i(t+b))|_{t=0}=i\in T_{(a+ib)}\C$. Now,
$${(\di\Psi_{(x,y)})}_{a+ib}(i)=\der{}{t}\Psi_{(x,y)}(a+i(t+b))\big|_{t=0}=(i_{N_1},i_{N_2})=(1_{N_2}, -1_{N_1})=J(1_{N_1},1_{N_2})\ .$$
The result follows observing that $i=i\cdot1$ is the action of the complex structure of $\C$ on the vector $1$.\\

\item Using the decomposition of $TN_i\simeq H_{i}\oplus \g{g}\simeq H_i\oplus \R^2$ and since $\Psi$ is holomorphic  we have $(\Psi_{a+ib})_*\circ J_i=J_i\circ (\Psi_{a+ib})_*$. For the vertical vectors we have, since they are invariant under $(\Psi_{a+ib})_*$ 
$$J((\Psi_{a+ib})_*\xi_{N_1},(\Psi_{a+ib})_*\eta_{N_2}))=((\Psi_{a+ib})_*\eta_{N_1},-(\Psi_{a+ib})_*\xi_{N_2})=(\eta_{N_1},-\xi_{N_2})\ ,$$
\newpage
\noindent and also,
$$(\Psi_{a+ib})_*J(\xi_{N_1},\eta_{N_2})=(\Psi_{a+ib})_*(\eta_{N_1},-\xi_{N_2})=(\eta_{N_1},-\xi_{N_2})\ .$$
\end{enumerate}
\end{proof}

\begin{rem}
    We observe that, in general, given $A\in GL(n,\Z)$, the action $\Psi:\C\times (N_1\times N_2)\to (N_1\times N_2)$ defined as
    $\Psi(a+ib,(x,y)):=(e^{2\pi i(A_{11}a+A_{21}b)}.x,e^{2\pi i(A_{12}a+A_{22}b)}.y)$ is holomorphic with respect to $J$, and induces a free holomorphic action of $\C/\Lambda_p$, where $\Lambda_p:=\{c+id\ :\ c,d\in\Z/p\}$, and $p=|\text{det}(A)|$.
\end{rem}

\begin{proof} of \textbf{Theorem \ref{thm2}}
    We find holomorphic charts of $N_1\times N_2$ and the holomorphic trivializations of $\pi:(N_1\times N_2)\to (P_1\times P_2)$ for the case of an $\sph^{1}$ action.

\subsection*{Holomorphic charts for $N_1\times N_2$} Let $U_1\subseteq P_1$ and $U_2\subseteq P_2$ open neighborhoods of $x_1\in U_1$ and $x_2\in U_2$ respectively. Then if $\pi_i:N_i\to P_i$, there are diffeomorphisms
$$\phi_i:\pi_i^{-1}(U_i)\to U_i\times \sph^1\ .$$
Let us denote with $\Sigma_i:=\phi^{-1}(U_i\times\{1\})$. Then we have 
$$T_{(x_1,x_2)}(N_1\times N_2)\simeq T_{(x_1,x_2)}(\Sigma_1\times \Sigma_2)\oplus T_{(x_1,x_2)}G.(x_1,x_2)\ .$$
Up to shrinking $\Sigma_i$ we can assume the previous relation to hold everywhere on $\Sigma_1\times \Sigma_2$.
Let $\Phi:\Sigma_1\times \Sigma_2\times G\to N_1\times N_2$ be the map
$$\Phi(s_1,s_2,g):=g.(s_1,s_2)\ .$$
This map is holomorphic thanks to proposition \ref{psihol1} and a local diffeomorphism.

\begin{lemma} The map $\Phi$ defined above is injective and hence a biholomorphism.
\end{lemma}

\begin{proof}
    Suppose by contradiction that $\Phi$ is not injective. Then pick $$\{(s_1^n,s_2^n)\}_{n\in \N}\subset \Sigma_1\times\Sigma_2\supset\{(\tilde s_1^n,\tilde s_2^n)\}_{n\in \N}\ ,$$ such that 
$(s_1^n,s_2^n)\neq (\tilde s_1^n,\tilde s_2^n)$, $g_n$ for all $n\in \N$ and 
    $$(s_1^n,s_2^n)\overset{n\to\infty}{\longrightarrow}(s_1,s_2)\overset{n\to\infty}{\longleftarrow}(\tilde s_1^n,\tilde s_2^n)\ .$$ 
    We can then find $\{g_n\}_{n\in\N}\subset G\supset \{\tilde g_n\}_{n\in\N}$ such that 
    $$g_n.(s_1^n,s_2^n)=\tilde g_n.(\tilde s_1^n,\tilde s_2^n)\ .$$
    Therefore if $\theta_n:=\tilde g_n^{-1}g_n$ we have $$\Phi(s_1^n,s_2^n,\theta_n)=\Phi(\tilde s_1^n,\tilde s_2^n, 1)=(\tilde s_1^n,\tilde s_2^n)\ .$$
    We conclude that $\Phi(s_1^n,s_2^n,\theta_n)\overset{n\to \infty}{\longrightarrow}(s_1,s_2)$ and by properness of the action $\theta_n\overset{n\to \infty}{\longrightarrow}1$. But this means that $\Phi$ is not injective in a neighborhood of $(s_1,s_2,1)$, wich is a contradiction.
\end{proof}

\subsection*{Transition functions} Let $\tilde U_i\subseteq P_i$ with $\tilde U_i\cap U_i\neq\emptyset$ and diffeomorphism 
$$\tilde\phi_i:\pi_i^{-1}(\tilde U_i)\to \tilde U_i\times \sph^1\ ,$$
and denote $\tilde \Sigma_i:=\tilde\phi_i^{-1}(\tilde U_i\times\{1\})$.
Again we have a biholomorphism with the image $\tilde\Phi:\tilde\Sigma_1\times\tilde\Sigma_2\times G\to N_1\times N_2$.
Let us use the following notation for simplicity
$\Sigma:=\Sigma_1\times\Sigma_2$, $\tilde\Sigma:=\tilde\Sigma_1\times\Sigma_2$, $U:=U_1\times U_2$, $\tilde U:=\tilde U_1\times \tilde U_2$,
$$\Phi^{-1}:\Phi(\Sigma\times G)\to \Sigma\times G,\ n\mapsto (\Pi^{-1}(\pi(n)),\Theta(n))\ ,$$
$$\tilde\Phi^{-1}:\tilde\Phi(\tilde\Sigma\times G)\to \tilde\Sigma\times G,\ n\mapsto (\tilde\Pi^{-1}(\pi(n)),\tilde\Theta(n))\ ,$$
Where $\Pi:=\pi|_{\Sigma}$ and $\tilde \Pi:=\pi|_{\tilde \Sigma}$ are smooth invertible maps. Then we can define trivializations $H:U\times G\to\pi^{-1}(U) $ as $$H(p,g):=\Phi(\Pi^{-1}(p),g)\ .$$
Similarly, $$\tilde H(p,g):=\tilde \Phi(\tilde\Pi^{-1}(p),g)\ .$$
We have $H^{-1}(n)=(\pi(n),\Theta (n))$. Now, 
$$(H^{-1}\circ \tilde H)(p,g)=H^{-1}(\tilde\Phi(\tilde \Pi^{-1}(p),g)=(\pi(\tilde\Phi(\tilde \Pi^{-1}(p),g)),\Theta(\tilde\Phi(\tilde \Pi^{-1}(p),g)))=(p,\Theta(\tilde\Phi(\tilde \Pi^{-1}(p),g)))\ .$$
Let $h(p)$ such that $\tilde \Phi(\tilde\Pi^{-1}(p),e)=h(p).\Phi(\Pi^{-1}(p),e)$, where $e=(1,1)$. Then, by equivariance, 
$$(H^{-1}\circ \tilde H)(p,g)=(p,\Theta(\tilde\Phi(\tilde \Pi^{-1}(p),g)))=(p,\Theta(g.\tilde\Phi (\Pi^{-1}(p),e)))=\qquad $$ $$=(p,\Theta( gh(p).\Phi( \Pi^{-1}(p),e)))=(p,\Theta(\Phi( \Pi^{-1}(p),gh(p))))=(p,gh(p))\  .$$

In summary, we have proved that the transition functions of the principal fiber bundle $\pi:(N_1\times N_2)\to (P_1\times P_2)$ are holomorphic. Hence $N_1\times N_2$ is a complex manifold, since the map $\Phi$ is holomorphic with respect to $J$. It also follows that the natural complex structure induced by the complex atlas is in fact $J$. Moreover, if we realize $P_1\times P_2$ as Meyer--Marsden--Weinstein reduction by an action of $\sph^{1}\times\overset{k-\text{times}}{\cdot\ \cdot\ \cdot}\times \sph^{1}$, then the same argument proves that $\pi:(N_1\times N_2)\to (P_1\times P_2)$ is a $\left(\mathbb{T}^{2}\times\overset{k-\text{times}}{\cdot\ \cdot\ \cdot}\times \mathbb{T}^{2}\right)$-principal bundle. The transition functions are then holomorphic functions and then, as before $N_1\times N_2$ admits a complex atlas and $J$ is precisely the complex structure induced by the atlas.
\end{proof}

\section{Applications}

\subsection{Complex Stiefel manifolds} Consider the action of the unitary group in complex dimension $k$, $U(k)$, on the space of linear maps $\text{Hom}(\C^k,\C^{n})$, for $n\geq k$. 
$$\Psi:U(k)\times \text{Hom}(\C^k,\C^{n})\to \text{Hom}(\C^k,\C^{n})\ ,$$
$$\Psi(u,A):=A\circ u^\dagger\ .$$
This is a proper left action\footnote{The action is proper since U(k) is compact, and is a left action since $\Psi(u_2u_1,A)=A\circ(u_2u_1)^\dagger=A\circ u_1^\dagger\circ u_2^\dagger=\Psi(u_2,\Psi(u_1,A))$.}. The action is also isometric with respect to the Hermitian metric metric on $\text{Hom}(\C^k,\C^n)$ given by $$h(A,B):=\text{Tr}(A\circ B^\dagger)\ .$$ 
The symplectic form associated, $\omega(A,B):=\im{\text{Tr}(A\circ B^\dagger)}$, is exact and equal to the differential of the Liouville 1-form $\theta|_A(B)=\im{\text{Tr}(A\circ B^\dagger)}$.
The infinitesimal generator of the action at $A\in\text{Hom}(\C^k,\C^n)$ for $\xi\in\g{u}(k)$ is given by $\der{}{t}\left(A\circ\exp{t\xi}\right)\big|_{t=0}=A\circ\xi$ and then we have
$$\dcp{\mu(A)}{\xi}=\theta_A(A\circ\xi)=\im{\text{Tr}(A\circ \xi^\dagger \circ A^\dagger)}=-\im{\text{Tr}(A^\dagger\circ A\circ \xi)}=i\text{Tr}(A^\dagger\circ A\circ\xi)\ .$$
By duality then, $$\dcp{\mu(A)}{\xi}=\text{Tr}(\mu(A)\circ\xi^\dagger)=-\text{Tr}(\mu(A)\circ\xi)=i\text{Tr}(A^\dagger\circ A\circ\xi)\ ,$$
and hence
$$\mu(A)=i(A^\dagger\circ A)\in \g{u}(k)^*\ .$$
This momentum map is $\Ad_{U(k)}$-equivariant, indeed $$\mu(A\circ u^\dagger)=i(u\circ A^\dagger\circ A\circ u^\dagger)=u\circ \mu(A)\circ u^\dagger\ .$$
Notice that the level set of $\mu$ relative to the element $i\mathbb{I}\in\g{u}(k)^*$ is $\mu^{-1}(i\mathbb{I})=\{A\in\text{Hom}(\C^k,\C^n)\ :\ A^\dagger\circ A=\mathbb{I}\}$ is equal to the set of unitary $k$-frames in $\C^n$, the so called Stiefel manifold $V_k(\C^n)$, and that the isotropy subgroup of any element, which is the whole $U(k)$, acts freely on it. We also know that the quotient space is given by the complex Grassmannian $\text{Gr}_\C(k,n)$, which is Kähler. Thanks to our construction, the product 
$V_k(\C^n)\times V_k(\C^n)$ carries a non integrable almost complex structure. We can also consider a $p$-dimensional torus in $U(k)$, $p\leq k$, and obtain a free action by $\C^p/\Lambda$ on $V_k(\C^n)\times V_k(\C^n)$, which is then foliated by tori. In particular, if $p=1$, the orbits are $J$-holomorphic curves (see \cite{McSal}).

\subsection{Infinite Calabi-Eckmann manifolds} Consider an Hilbert space $\mathcal{H}$ and its unit sphere $S(\mathcal{H})$. Theorem \ref{thm2} tells us that the product $\sph^{2n+1}\times S(\mathcal{H})$ is holomorphic with complex structure $J$ and then we have the following.

\begin{thm}
     There cannot be a $J$-invarinat symplectic structure $\omega$ on $\sph^{2n+1}\times S(\mathcal{H})$ such that $\omega(\cdot,J\cdot)$ is Riemannian. In other words $\sph^{2n+1}\times S(\mathcal{H})$ cannot be Kähler with respect to the complex structure $J$.
\end{thm} 

\begin{proof}
    If such a symplectic form $\omega$ does exist, the complex submanifolds $j:\sph^{2n+1}\times \sph^{2m+1}\hookrightarrow \sph^{2n+1}\times S(\mathcal{H})$ wolud be Kähler with respect to $j^*\omega$. This is a contradiction since $\sph^{2n+1}\times \sph^{2m+1}$ cannot be Kähler, e.g. as we have already recalled, for cohomological reasons, when $n\geq 1 $ or $m\geq 1$.
    
\end{proof}

\newpage

\end{document}